\documentclass{article}
\usepackage{amssymb,amsmath,theorem,euscript}

\input{epsf}

\newcounter{sec}

\def\sm{\smallskip}

\newcounter{punct}[sec]

\def\punct{\refstepcounter{punct}{\arabic{sec}.\arabic{punct}.  }}

\def\COUNTERS{\addtocounter{sec}{1}
              \setcounter{punct}{0}
          \setcounter{equation}{0}
          \setcounter{theorem}{0}
                  }
\newtheorem{theorem}{Theorem}[sec]

\newtheorem{lemma}[theorem]{Lemma}
\newtheorem{corollary}[theorem]{Corollary}

\begin{document}

 \def\ov{\overline}
\def\wt{\widetilde}
 \newcommand{\rk}{\mathop {\mathrm {rk}}\nolimits}
\newcommand{\Aut}{\mathop {\mathrm {Aut}}\nolimits}
\newcommand{\Out}{\mathop {\mathrm {Out}}\nolimits}
\renewcommand{\Re}{\mathop {\mathrm {Re}}\nolimits}
\renewcommand{\Im}{\mathop {\mathrm {Im}}\nolimits}
 \newcommand{\tr}{\mathop {\mathrm {tr}}\nolimits}
  \newcommand{\Hom}{\mathop {\mathrm {Hom}}\nolimits}
   \newcommand{\diag}{\mathop {\mathrm {diag}}\nolimits}
   \newcommand{\supp}{\mathop {\mathrm {supp}}\nolimits}
 \newcommand{\im}{\mathop {\mathrm {im}}\nolimits}
 \newcommand{\Res}{\mathop {\mathrm {Res}}\nolimits}   
   
\def\Br{\mathrm {Br}}
\def\Fl{\mathrm {Fl}}
 %%%1. ClASSICAL GROUPS
\def\SL{\mathrm {SL}}
\def\SU{\mathrm {SU}}
\def\GL{\mathrm {GL}}
\def\U{\mathrm U}
\def\OO{\mathrm O}
 \def\Sp{\mathrm {Sp}}
 \def\SO{\mathrm {SO}}
\def\SOS{\mathrm {SO}^*}
 \def\Diff{\mathrm{Diff}}
 \def\Vect{\mathfrak{Vect}}
\def\PGL{\mathrm {PGL}}
\def\PU{\mathrm {PU}}
\def\PSL{\mathrm {PSL}}
\def\Symp{\mathrm{Symp}}
\def\End{\mathrm{End}}
\def\Mor{\mathrm{Mor}}
\def\Aut{\mathrm{Aut}}
 \def\PB{\mathrm{PB}}
 \def\cA{\mathcal A}
\def\cB{\mathcal B}
\def\cC{\mathcal C}
\def\cD{\mathcal D}
\def\cE{\mathcal E}
\def\cF{\mathcal F}
\def\cG{\mathcal G}
\def\cH{\mathcal H}
\def\cJ{\mathcal J}
\def\cI{\mathcal I}
\def\cK{\mathcal K}
 \def\cL{\mathcal L}
\def\cM{\mathcal M}
\def\cN{\mathcal N}
 \def\cO{\mathcal O}
\def\cP{\mathcal P}
\def\cQ{\mathcal Q}
\def\cR{\mathcal R}
\def\cS{\mathcal S}
\def\cT{\mathcal T}
\def\cU{\mathcal U}
\def\cV{\mathcal V}
 \def\cW{\mathcal W}
\def\cX{\mathcal X}
 \def\cY{\mathcal Y}
 \def\cZ{\mathcal Z}
%%% END MATHCAL %%%%%%%%%%%%%%%%%%%%%%%%%%%%%%%%% %%%%%%%%%%%%%%%%%%%%%%%%%%%%%%%% %%%
\def\0{{\ov 0}}
 \def\1{{\ov 1}}
 %%%%%%%%%%%%%%%%%%%%%%%%%%%% %%%%%%%%%%%%%%%%%%%%%%%%%%%%%%%%%%% %%% BEGIN GOTIC
 \def\frA{\mathfrak A}
 \def\frB{\mathfrak B}
\def\frC{\mathfrak C}
\def\frD{\mathfrak D}
\def\frE{\mathfrak E}
\def\frF{\mathfrak F}
\def\frG{\mathfrak G}
\def\frH{\mathfrak H}
\def\frI{\mathfrak I}
 \def\frJ{\mathfrak J}
 \def\frK{\mathfrak K}
 \def\frL{\mathfrak L}
\def\frM{\mathfrak M}
 \def\frN{\mathfrak N} \def\frO{\mathfrak O} \def\frP{\mathfrak P} \def\frQ{\mathfrak Q} \def\frR{\mathfrak R}
 \def\frS{\mathfrak S} \def\frT{\mathfrak T} \def\frU{\mathfrak U} \def\frV{\mathfrak V} \def\frW{\mathfrak W}
 \def\frX{\mathfrak X} \def\frY{\mathfrak Y} \def\frZ{\mathfrak Z} \def\fra{\mathfrak a} \def\frb{\mathfrak b}
 \def\frc{\mathfrak c} \def\frd{\mathfrak d} \def\fre{\mathfrak e} \def\frf{\mathfrak f} \def\frg{\mathfrak g}
 \def\frh{\mathfrak h} \def\fri{\mathfrak i} \def\frj{\mathfrak j} \def\frk{\mathfrak k} \def\frl{\mathfrak l}
 \def\frm{\mathfrak m} \def\frn{\mathfrak n} \def\fro{\mathfrak o} \def\frp{\mathfrak p} \def\frq{\mathfrak q}
 \def\frr{\mathfrak r} \def\frs{\mathfrak s} \def\frt{\mathfrak t} \def\fru{\mathfrak u} \def\frv{\mathfrak v}
 \def\frw{\mathfrak w} \def\frx{\mathfrak x} \def\fry{\mathfrak y} \def\frz{\mathfrak z} \def\frsp{\mathfrak{sp}}
 %% This is Lie algebra %%% END GOTIC
%%%%%%%%%%%%%%%%%%%%%%%%%%%%%%%% %%%%%%%%%%%%%%%%%%%%%%%%%%%%%%%%%
%%% BEGIN MATHBF
 \def\bfa{\mathbf a} \def\bfb{\mathbf b} \def\bfc{\mathbf c} \def\bfd{\mathbf d} \def\bfe{\mathbf e} \def\bff{\mathbf f}
 \def\bfg{\mathbf g} \def\bfh{\mathbf h} \def\bfi{\mathbf i} \def\bfj{\mathbf j} \def\bfk{\mathbf k} \def\bfl{\mathbf l}
 \def\bfm{\mathbf m} \def\bfn{\mathbf n} \def\bfo{\mathbf o} \def\bfp{\mathbf p} \def\bfq{\mathbf q} \def\bfr{\mathbf r}
 \def\bfs{\mathbf s} \def\bft{\mathbf t} \def\bfu{\mathbf u} \def\bfv{\mathbf v} \def\bfw{\mathbf w} \def\bfx{\mathbf x}
 \def\bfy{\mathbf y} \def\bfz{\mathbf z} \def\bfA{\mathbf A} \def\bfB{\mathbf B} \def\bfC{\mathbf C} \def\bfD{\mathbf D}
 \def\bfE{\mathbf E} \def\bfF{\mathbf F} \def\bfG{\mathbf G} \def\bfH{\mathbf H} \def\bfI{\mathbf I} \def\bfJ{\mathbf J}
 \def\bfK{\mathbf K} \def\bfL{\mathbf L} \def\bfM{\mathbf M} \def\bfN{\mathbf N} \def\bfO{\mathbf O} \def\bfP{\mathbf P}
 \def\bfQ{\mathbf Q} \def\bfR{\mathbf R} \def\bfS{\mathbf S} \def\bfT{\mathbf T} \def\bfU{\mathbf U} \def\bfV{\mathbf V}
 \def\bfW{\mathbf W} \def\bfX{\mathbf X} \def\bfY{\mathbf Y} \def\bfZ{\mathbf Z} \def\bfw{\mathbf w}
 %%% END MATHBF
%%%%%%%%%%%%%%%%%%%%%%%%%%%%%%% %%%%%%%%%%%%%%%%%%%%%%%%%%%%%%%%%
 %%% BEGIN MATHBB
 \def\R {{\mathbb R }} \def\C {{\mathbb C }} \def\Z{{\mathbb Z}} \def\H{{\mathbb H}} \def\K{{\mathbb K}}
 \def\N{{\mathbb N}} \def\Q{{\mathbb Q}} \def\A{{\mathbb A}} \def\T{\mathbb T} \def\P{\mathbb P} \def\G{\mathbb G}
 \def\bbA{\mathbb A} \def\bbB{\mathbb B} \def\bbD{\mathbb D} \def\bbE{\mathbb E} \def\bbF{\mathbb F} \def\bbG{\mathbb G}
 \def\bbI{\mathbb I} \def\bbJ{\mathbb J} \def\bbL{\mathbb L} \def\bbM{\mathbb M} \def\bbN{\mathbb N} \def\bbO{\mathbb O}
 \def\bbP{\mathbb P} \def\bbQ{\mathbb Q} \def\bbS{\mathbb S} \def\bbT{\mathbb T} \def\bbU{\mathbb U} \def\bbV{\mathbb V}
 \def\bbW{\mathbb W} \def\bbX{\mathbb X} \def\bbY{\mathbb Y} \def\kappa{\varkappa} \def\epsilon{\varepsilon}
 \def\phi{\varphi} \def\le{\leqslant} \def\ge{\geqslant}

\def\UU{\bbU}
\def\Mat{\mathrm{Mat}}
\def\Herm{\mathrm{Herm}}
\def\AHerm{\mathrm{AHerm}}
\def\tto{\rightrightarrows}

\def\Gms{\mathrm {Gms}}
\def\Ams{\mathrm {Ams}}
\def\Isom{\mathrm {Isom}}

\def\Gr{\mathrm{Gr}}

\def\graph{\mathrm{graph}}

\def\O{\mathrm{O}}

\def\la{\langle}
\def\ra{\rangle}

%\begin{document}

 \def\ov{\overline}
\def\wt{\widetilde}

\renewcommand{\Re}{\mathop {\mathrm {Re}}\nolimits}
\def\Br{\mathrm {Br}}

 %%%1. ClASSICAL GROUPS
 \def\Isom{\mathrm {Isom}}
\def\SL{\mathrm {SL}}
\def\SU{\mathrm {SU}}
\def\GL{\mathrm {GL}}
\def\U{\mathrm U}
\def\OO{\mathrm O}
 \def\Sp{\mathrm {Sp}}
 \def\SO{\mathrm {SO}}
\def\SOS{\mathrm {SO}^*}
 \def\Diff{\mathrm{Diff}}
 \def\Vect{\mathfrak{Vect}}
\def\PGL{\mathrm {PGL}}
\def\PU{\mathrm {PU}}
\def\PSL{\mathrm {PSL}}
\def\Symp{\mathrm{Symp}}
\def\End{\mathrm{End}}
\def\Mor{\mathrm{Mor}}
\def\Aut{\mathrm{Aut}}
 \def\PB{\mathrm{PB}}
 \def\cA{\mathcal A}
\def\cB{\mathcal B}
\def\cC{\mathcal C}
\def\cD{\mathcal D}
\def\cE{\mathcal E}
\def\cF{\mathcal F}
\def\cG{\mathcal G}
\def\cH{\mathcal H}
\def\cJ{\mathcal J}
\def\cI{\mathcal I}
\def\cK{\mathcal K}
 \def\cL{\mathcal L}
\def\cM{\mathcal M}
\def\cN{\mathcal N}
 \def\cO{\mathcal O}
\def\cP{\mathcal P}
\def\cQ{\mathcal Q}
\def\cR{\mathcal R}
\def\cS{\mathcal S}
\def\cT{\mathcal T}
\def\cU{\mathcal U}
\def\cV{\mathcal V}
 \def\cW{\mathcal W}
\def\cX{\mathcal X}
 \def\cY{\mathcal Y}
 \def\cZ{\mathcal Z}
%%% END MATHCAL %%%%%%%%%%%%%%%%%%%%%%%%%%%%%%%%% %%%%%%%%%%%%%%%%%%%%%%%%%%%%%%%% %%%
\def\0{{\ov 0}}
 \def\1{{\ov 1}}
 %%%%%%%%%%%%%%%%%%%%%%%%%%%% %%%%%%%%%%%%%%%%%%%%%%%%%%%%%%%%%%% %%% BEGIN GOTIC
 \def\frA{\mathfrak A}
 \def\frB{\mathfrak B}
\def\frC{\mathfrak C}
\def\frD{\mathfrak D}
\def\frE{\mathfrak E}
\def\frF{\mathfrak F}
\def\frG{\mathfrak G}
\def\frH{\mathfrak H}
\def\frI{\mathfrak I}
 \def\frJ{\mathfrak J}
 \def\frK{\mathfrak K}
 \def\frL{\mathfrak L}
\def\frM{\mathfrak M}
 \def\frN{\mathfrak N} \def\frO{\mathfrak O} \def\frP{\mathfrak P} \def\frQ{\mathfrak Q} \def\frR{\mathfrak R}
 \def\frS{\mathfrak S} \def\frT{\mathfrak T} \def\frU{\mathfrak U} \def\frV{\mathfrak V} \def\frW{\mathfrak W}
 \def\frX{\mathfrak X} \def\frY{\mathfrak Y} \def\frZ{\mathfrak Z} \def\fra{\mathfrak a} \def\frb{\mathfrak b}
 \def\frc{\mathfrak c} \def\frd{\mathfrak d} \def\fre{\mathfrak e} \def\frf{\mathfrak f} \def\frg{\mathfrak g}
 \def\frh{\mathfrak h} \def\fri{\mathfrak i} \def\frj{\mathfrak j} \def\frk{\mathfrak k} \def\frl{\mathfrak l}
 \def\frm{\mathfrak m} \def\frn{\mathfrak n} \def\fro{\mathfrak o} \def\frp{\mathfrak p} \def\frq{\mathfrak q}
 \def\frr{\mathfrak r} \def\frs{\mathfrak s} \def\frt{\mathfrak t} \def\fru{\mathfrak u} \def\frv{\mathfrak v}
 \def\frw{\mathfrak w} \def\frx{\mathfrak x} \def\fry{\mathfrak y} \def\frz{\mathfrak z} \def\frsp{\mathfrak{sp}}
 %% This is Lie algebra %%% END GOTIC
%%%%%%%%%%%%%%%%%%%%%%%%%%%%%%%% %%%%%%%%%%%%%%%%%%%%%%%%%%%%%%%%%
%%% BEGIN MATHBF
 \def\bfa{\mathbf a} \def\bfb{\mathbf b} \def\bfc{\mathbf c} \def\bfd{\mathbf d} \def\bfe{\mathbf e} \def\bff{\mathbf f}
 \def\bfg{\mathbf g} \def\bfh{\mathbf h} \def\bfi{\mathbf i} \def\bfj{\mathbf j} \def\bfk{\mathbf k} \def\bfl{\mathbf l}
 \def\bfm{\mathbf m} \def\bfn{\mathbf n} \def\bfo{\mathbf o} \def\bfp{\mathbf p} \def\bfq{\mathbf q} \def\bfr{\mathbf r}
 \def\bfs{\mathbf s} \def\bft{\mathbf t} \def\bfu{\mathbf u} \def\bfv{\mathbf v} \def\bfw{\mathbf w} \def\bfx{\mathbf x}
 \def\bfy{\mathbf y} \def\bfz{\mathbf z} \def\bfA{\mathbf A} \def\bfB{\mathbf B} \def\bfC{\mathbf C} \def\bfD{\mathbf D}
 \def\bfE{\mathbf E} \def\bfF{\mathbf F} \def\bfG{\mathbf G} \def\bfH{\mathbf H} \def\bfI{\mathbf I} \def\bfJ{\mathbf J}
 \def\bfK{\mathbf K} \def\bfL{\mathbf L} \def\bfM{\mathbf M} \def\bfN{\mathbf N} \def\bfO{\mathbf O} \def\bfP{\mathbf P}
 \def\bfQ{\mathbf Q} \def\bfR{\mathbf R} \def\bfS{\mathbf S} \def\bfT{\mathbf T} \def\bfU{\mathbf U} \def\bfV{\mathbf V}
 \def\bfW{\mathbf W} \def\bfX{\mathbf X} \def\bfY{\mathbf Y} \def\bfZ{\mathbf Z} \def\bfw{\mathbf w}
 %%% END MATHBF
%%%%%%%%%%%%%%%%%%%%%%%%%%%%%%% %%%%%%%%%%%%%%%%%%%%%%%%%%%%%%%%%
 %%% BEGIN MATHBB
 \def\R {{\mathbb R }} \def\C {{\mathbb C }} \def\Z{{\mathbb Z}} \def\H{{\mathbb H}} \def\K{{\mathbb K}}
 \def\N{{\mathbb N}} \def\Q{{\mathbb Q}} \def\A{{\mathbb A}} \def\T{\mathbb T} \def\P{\mathbb P} \def\G{\mathbb G}
 \def\bbA{\mathbb A} \def\bbB{\mathbb B} \def\bbD{\mathbb D} \def\bbE{\mathbb E} \def\bbF{\mathbb F} \def\bbG{\mathbb G}
 \def\bbI{\mathbb I} \def\bbJ{\mathbb J} \def\bbL{\mathbb L} \def\bbM{\mathbb M} \def\bbN{\mathbb N} \def\bbO{\mathbb O}
 \def\bbP{\mathbb P} \def\bbQ{\mathbb Q} \def\bbS{\mathbb S} \def\bbT{\mathbb T} \def\bbU{\mathbb U} \def\bbV{\mathbb V}
 \def\bbW{\mathbb W} \def\bbX{\mathbb X} \def\bbY{\mathbb Y} \def\kappa{\varkappa} \def\epsilon{\varepsilon}
 \def\phi{\varphi} \def\le{\leqslant} \def\ge{\geqslant}

\def\UU{\bbU}
\def\Mat{\mathrm{Mat}}
\def\tto{\rightrightarrows}

\def\Gr{\mathrm{Gr}}

\def\graph{\mathrm{graph}}

\def\O{\mathrm{O}}

\def\la{\langle}
\def\ra{\rangle}

\begin{center}
 \Large \bf
 Radial parts of Haar measures 
\\
 and probability distributions on the space
\\
 of rational matrix-valued functions
 
 \large\sc
 
 \bigskip
 Yu.A.Neretin%
\footnote{The research was carried out at the IITP RAS at the expense of the Russian Foundation for Sciences 
(project № 14-50-00150).}
 
\end{center}

{\small Consider the space  $\cC$ of conjugacy classes of a unitary group
$\U(n+m)$ with respect to a smaller unitary group   $\U(m)$. 
It is known that for any element 
of the space  $\cC$ we can assign canonically a matrix-valued rational function
on the Riemann sphere (a Livshits characteristic function).
In the paper we write an explicit expression for the natural measure 
on $\cC$ obtained as the pushforward of the Haar measure of the group 
$\U(n+m)$ in the terms of characteristic functions.}

\section{The statement}

\COUNTERS

{\bf\punct The purpose of the paper.}
There is a wide literature 
(see, e.g.,  \cite{PV}, \cite{Sod}, \cite{ST},  \cite{HKM}
and further references in these works)
on Gaussian random functions, the topic arises at least to the Payley--Wiener book
 \cite{PW}, Chapter 10. Relatively recently M.~Krishnapur 
\cite{Kri} started investigation of random matrix-valued holomorphic functions.

In the present paper we consider measures on the space of rational matrix-valued functions
on the Riemann sphere.
The origin of the question under a discussion is the following.
To be definite, consider an unitary group
$\U(n)$. The distribution of eigenvalues of unitary matrices is a measure on the set of 
$n$-point subsets on the circle with a density of the form
 $C\cdot \prod_{k<l} |z_k-z_l|^2$
(The Hermann Weyl formula, see, e.g., \cite{Hua}, formula (3.2.2), or \cite{Far}, Theorem 11.2.1).
There exists a zoo of similar formulas, a usual corresponding term is 
'radial parts of Haar measures on symmetric spaces'.
Namely, we consider a Riemannian symmetric space 
$G/K$, the space of double cosets  $K\setminus G/K$,
and the pushforward of the Haar measure under the map  $G\to K\setminus G/K$ 
(in the example with unitary group, we have  $G=\U(n)\times \U(n)$,
and $K$ is the diagonal subgroup $\U(n)$),
general formulas are contained in 
 \cite{Hel}, Propositions   X.1.17, X.1.19. 

 This important topic has numerous applications and continuations.
 However, its extensions to other pairs of groups and subgroups are almost absent%
\footnote{The case when $G=\SU(2)\times \dots \times \SU(2)$, and   $K=\SU(2)$ is the diagonal subgroup,
is examined in \cite{Ner-SU2}.}
$G\supset K$. 
One of obstacles for such extensions are difficulties related to descriptions
of double coset and conjugacy classes spaces.

In one case such description is known during a long time.
Namely, for conjugacy classes of a unitary group 
 $\U(n+m)$ with respect to a smaller subgroup
 $\U(m)$ the solution is given in the terms of  Livshits characteristic functions
 (see below). In
\cite{Ner-book} there was proposed a counterpart of characteristic functions for double cosets 
of $\U(n+m)$ by $\O(m)$, see  also \cite{Olsh}. In \cite{Ner-char},  \cite{Ner-pp} characteristic functions were
constructed for a wide class  pairs group--subgroup. 
Since characteristic functions from  \cite{Ner-char} were
originated from representation theory, see 
\cite{Ner-faa}, the question about radial parts of Haar measures
arises naturally. 

In the present paper we get an explicit formula for radial part of Haar measure
in the case of conjugacy classes of
 $\U(n+m)$ by $\U(m)$. 

\sm

{\bf\punct Livshits characteristic function.%
\label{ss:livshits}}
Let $g=\begin{pmatrix}
          \alpha&\beta\\ \gamma&\delta
         \end{pmatrix}$
         be a block matrix of size
  $(n+m)\times(n+m)$. The {\it characteristic function}
      of the matrix  $g$
         (another  term is a  {\it transfer-function})
         is a function on  $\C$ defined by the formula  
         \begin{equation}
         \chi(\lambda):=\alpha+\lambda\beta(1-\lambda\delta)^{-1}\gamma.
         \label{eq:livshits}
         \end{equation}
         This function takes values in the space of matrices of size
 $n\times n$.
 It is easy to see that the function does not changes under conjugations 
 of the matrix 
 $g$ by matrices of the form 
$\begin{pmatrix} 1&0\\0& u  \end{pmatrix}$. I other words, $\chi$ is an invariant
of conjugacy classes of 
$\U(n+m)$ by $\U(m)$.

It is known that any rational matrix-valued function that has not a singularity
at 
$\lambda=0$, is a characteristic function of a matrix.
On a reconstruction  of the  matrix (more precisely, of the conjugacy class) $g$ from its characteristic function, see, e.g.,
\cite{Dym}, Chapter 19.
This reconstruction is not unique%
\footnote{However a matrix of a minimal possible size 
with a given characteristic function is 
unique up to a conjugation, see the textbook Dym  \cite{Dym}, Chapter 19.
Moreover an element of a categorical quotient of
 $\GL(n+m,\C)$ by $\GL(m,\C)$ can be reconstructed from the characteristic function
 in a unique way
(see \cite{Ner-p}).}.

\sm

Now let a matrix 
 $g$ be unitary. Then the characteristic function satisfies
 the following properties:
\sm

1) For $|\lambda|=1$ values of  $\chi(\lambda)$ are unitary matrices%
\footnote{A proof of this and the next statement are contained in 
Remark in .\ref{ss:lemma}.}.

\sm

2) For  $|\lambda|<1$ we have%
\footnote{Thus characteristic functions are matrix-valued analogs
of  {\it interior functions}, which are well-known in classical theory
of analytic functions, see, e.g.,  \cite{Gar}.
The theory of matrix-valued interior functions was developed in 
\cite{Pot}. On interior functions of matrix argument,
see  \cite{Ner-char}.}
$\|\chi(\lambda)\|\le 1$; by the Riemann--Schwarz reflection principle,
for $|\lambda|>1$
we have $\|\chi(\lambda)^{-1}\|\le 1$.

\sm

3) Represent $\det \chi(\lambda)$ as an irreducible fraction, 
$\frac{u(\lambda)}{v(\lambda)}$. It is easy to show that 
degrees of polynomials  $u(\lambda)$,  $v(\lambda)$ are $\le m$.
Indeed this determinant can be written in the form%
\footnote{For this, it is sufficient to apply 
the usual formula for a determinant of a block matrix,
see below (\ref{eq:block-determinant}).}
\begin{equation}
\det\chi(\lambda)=
\frac{\det\begin{pmatrix}
           \alpha&-\lambda \beta\\
           \gamma&1-\lambda\delta
          \end{pmatrix}
}{\det (1-\lambda\delta) }
.
\label{eq:det}
\end{equation}

\sm

Conversely, any rational function on $\C$
taking values in the space of 
$n\times n$-matrices, satisfying the above-listed properties, is a characteristic function 
of a certain unitary matrix of size  
 $n+m$. {\it Denote the set of such functions by  $\cR_{n}(m)$.}

\sm

Next, consider block matrices of size
 $n+m=n+(m_1+m_2)$ having the structure
\begin{equation}
g=\begin{pmatrix}
 \alpha&\beta&0\\
 \gamma&\delta&0\\
 0&0&\kappa
 \end{pmatrix}.
 \label{eq:split}
\end{equation}
It is easy to see that a characteristic function does not depend on the block
$\kappa$, on the other hand eigenvalues of the matrix 
$ \kappa$ are invariants of a conjugacy class.
It is easy to show that any matrix  $g\in\U(n+m)$
by conjugations by elements of
 $\U(m)$ can be reduced to the form  (\ref{eq:split}), where $\|\delta\|<1$.

\sm

Finally,
{\it spectral data determining a conjugacy class of  $\U(n+m)$ 
by $\U(m)$ looks as follows: the characteristic function
and  the set of eigenvalues of the matrix
 $\delta$ that are contained in the circle  $|\lambda|=1$.}

\sm

{\sc Remark.} I do not know a source of literature
containing simple proofs of statements about {\it unitary}
matrices in a (very small) degree of generality that is necessary for us.
For arbitrary matrices a short exposition is contained in 
 \cite{Dym}, Chapter  19, see also
 \cite{Bart}. The unitary case requires additional arguments.
 The statements that are necessary for us is an extremely particular case of
 \cite{Bro}, Theorem 5.1. A restoring of a unitary matrix from a characteristic
 function also can be done by the Potapov method
 \cite{Pot}, using expansion of a rational matrix-valued inner function
 in a Blaschke product.
 \hfill $\boxtimes$
 
 \sm
 
 Anyway, for elements of
$\U(n+m)$ in a general position, the block  $\delta$ 
is purely contractive and conjugacy classes are determined by their
characteristic functions.
Thus we have the map  
 $$
 \U(n+m)\to \cR_n(m)
 ,$$
 and we wish to evaluate the image of the Haar measure under this map.
 The answer will be presented in terms of the set of $\lambda$, where 
 $\chi(\lambda)$ has an eigenvalue 
 $-1$, and the corresponding eigenvectors. 

\sm

{\bf \punct The density of the measure.} Consider an element of
the space $\cR_n(m)$ in a general position
(i.e., an element of a set of full measure, which will be detaches step-by-step).
Consider the set of all points 
 $t_k\in \C$ such that a matrix 
 $\chi(t_k)$ has an eigenvalue  $-1$. This set is contained in the circle %
\footnote{For elements of  $\U(n+m)$ in a general position we have  $\|\alpha\|<1$. But
$\alpha=\chi(0)$. On the other hand, on the unit circle $\|\chi(\lambda)\|=1$.
Applying the maximum principle to linear functionals on the space
of matrices we  easily get
$\|\chi(\lambda)\|<1$ inside the circle.}
 $|\lambda|=1$. In a general position the number of such points is  $m$.
Let us order the points  $t_k$ in some way, for instance, let us assume that 
$$
2\pi>\arg t_1>\dots > \arg  t_n>0.
$$
In a general position there is a unique corresponding eigenvector 
$c_k$,
$$
\chi(t_k) c_k=-c_k,\qquad c_k=(c_k^1,\dots,c_k^n)\in\C^n
.$$
To fix coordinates on
 $\cR_n(m)$, we normalize 
$c_k$ from the condition  
$$
\la \chi'(t_k)\, c_k,c_k  \ra=-t_k^{-1}.
$$
Next, we fix a phase of each vector 
$c_k$ by the assumption $c_k^1\ge 0$.

\begin{lemma}
\label{l:1}
The properties mentioned above really hold in a general position. 
 In a general position, the numbers $t_k$, the vectors $c_k$, and the unitary matrix 
 $$
 U:=\chi(-1)
 $$
 uniquely determine a characteristic function.
\end{lemma}

Compose a matrix
 $C$ of size $n\times m$ of vector-columns $c_k$.
Compose a diagonal  matrix $T$ of numbers $t_j$.

\begin{theorem}
\label{th:main}
The image of probabilistic Haar measure under the map
 $$\U(n+m)\to \cR_n(m)$$
 has the form
 \begin{multline}
 \theta_{n,m}\cdot
  \bigl|\det(1+T+C^*(1+U) C)\bigr|^{-2n-2m} \bigl|\det(1+U)\bigr|^{2m}\times\\\times
  \prod_{k=1}^m |1+t_k|^{2m+2n}
  \prod_{1\le k<l\le n}|t_k-t_l|^2
  \times\\\times
  d\sigma_n(U)\cdot \prod_{k=1}^m c_k^1\,dc_k^1 \cdot  \prod_{k=1}^m\prod_{j=2}^{n} d\Re c_k^j\,d\Im c_k^j\cdot 
  \prod_{k=1}^n \frac {dt_k}{it_k},
  \label{eq:main}
 \end{multline}
where $d\sigma_n(U)$ is the Haar measure on $\U(n)$ and
\begin{equation}
\theta_{n,m}=2^{-m}\pi^{-mn}\frac{\prod_{j=1}^{m+n-1} j!}{\prod_{j=1}^{n-1} j!\cdot\prod_{j=1}^{m} j!}
.
 \label{eq:constant}
\end{equation}
\end{theorem}

In view of this statement, we mention two works on random analytic functions.

First, measures on the space of scalar meromorphic functions in terms of
distributions of poles and residues were considered by Wigner
 \cite{Wig}, 1952. This is similar to the language of our Theorem \ref{th:main},
 moreover during a calculation in Subsect.\ref{ss:poles}
 we literally use distributions of poles and residues
(the main difference is the following:
in Wigner paper distributions
of poles and residues are determined by independent random variables,
in our case this is not so).
The class of meromorphic functions, which was considered in \cite{Wig}
can be send to the class of of (scalar-valued) interior functions
by a simple transformation.

Second, a measure on the space of interior functions was considered by Katsnelson
in
\cite{Kat}, 2002, it seems (at least from the first sight) that his object
is not similar to our considerations.

\sm

I am grateful to M.~Sodin, N.~Makarov, S.~Sodin for discussions of this topic.

\section{Proof}

\COUNTERS

Notation:

\sm

$\Mat_{k,l}$ is the space of complex matrices of size   $k\times l$;

\sm

$\Herm_k$ is the space of Hermitian matrices of size  $k$;

\sm

$\AHerm_k$ is the space of anti-Hermitian matrices  ($X^*=-X$) of size $k$;

\sm

$\T^k$ is the torus, i.e. the product of  $k$ circles with the corresponding group
structure;

\sm

$\R_+$ is the positive semi-axis.

\sm

{\bf\punct A preparation lemma.%
\label{ss:lemma}}
Let $T$ be a square matrix. Define its {\it Cayley transform}
by the formula
$$
g\mapsto H= -1+2(1+g)^{-1} =(1+g)^{-1}(1-g).
$$
This transform is inverse to itself and send unitary matrices to anti-Hermitian matrices.

\begin{lemma}
\label{l:2}
Consider a block matrix 
$g$ of size $n+m$  and perform the following chain 
of manipulations:
 
 \sm
 
{\rm 1)} we apply the Cayley transform and get a certain
matrix  $H$;

\sm
 
{\rm 2)} we take the characteristic function   $\phi(s)$ of the matrix $H$;

\sm
 
{\rm 3)} we apply the Cayley transform to  $\phi(s)$ and get
a function $\psi(s)$ taking values in  $\Mat_{n,n}$;
 
 \sm
 
{\rm 4)} we consider the function $\phi\Bigl(\frac{t+1}{t-1}\Bigr)$.

\sm

Then the resulting function coincides with the characteristic function
 $\chi_g(t)$ of the matrix $g$.
\end{lemma}

First, we recall that definitions of the Cayley transform and the characteristic function
can be formulated as follows.

\sm

 {$\bullet$ \it Let $g$ be a square matrix of size  $n$. Let vectors $u$, $v\in \C^n$ 
 be connected by the relation
 \begin{equation}
 (u-v)=g(u+v).
 \label{eq:guv}
 \end{equation}
 Then $v= Hu$, where $H$ is the Cayley transform of the matrix  $g$.}
 
 \sm
 
 {$\bullet$ \it Let $g$ be a block matrix of size  $n+m$. Fix $\lambda\in\C$. 
 Consider the set  $L$ of all pairs   $q$, $p\in \C^n$,
 for which there exists
 $x\in \C^m$ such that 
 \begin{equation}
 \begin{pmatrix}
  q\\x
 \end{pmatrix}=
 \begin{pmatrix}
  \alpha&\beta\\\gamma&\delta
 \end{pmatrix}
 \begin{pmatrix}
  p\\\lambda x
 \end{pmatrix}.
 \label{eq:qp}
 \end{equation}
 The subset 
 $L\subset \C^n\oplus \C^n$ obtained in such way is determined by the equation 
$$
q=\chi_g(\lambda)\, p
$$
for all
$\lambda$ except poles of the characteristic function.
 }
 
 \sm
 
 The first statement is evident, we simply write
(\ref{eq:guv}) in the form
$$
(1-g)u=(1+g)v.
$$
To verify the second statement, we write 
 (\ref{eq:qp}) as
 $$
 q=\alpha p+\lambda \beta x,\qquad x=\gamma p+\lambda\delta x.
 $$
 Eliminating
  $x$ we get the desired statement.
 
 \sm
 
 {\sc  Remark.} Now we can easily derive  properties  1)--2) of characteristic functions
 from Subsect.  \ref{ss:livshits}. Indeed, if a matrix  $g$ is unitary, then 
 $$
 \|q\|^2+\|x\|^2=\|p\|^2+\|\lambda x\|^2.
 $$
 Assuming $|\lambda|=1$, we get  $ \|q\|^2= \|p\|^2$, 
 this implies that the matrix
 $\chi(\lambda)$ is unitary. If $|\lambda|<1$,
 then $ \|q\|^2\ge \|p\|^2$, i.e., $\|\chi(\lambda)\|\le 1$. 
  	It is worth noting, that a proof using directly formula  
 (\ref{eq:livshits}) is not so easy. 
 \hfill
 $\boxtimes$
 
 \sm
 
 {\sc  Proof of the lemma.} Let $g$ be a block square matrix of size  $n+m$.
 Then
$v=Hu$ is equivalent to the equality 
 $$
 u-v=g(u+v).
 $$
 Next, $q=\phi(s)\,p$, if there exists a vector  $x$ such that
 $$
 \begin{pmatrix}
  p\\ s x
 \end{pmatrix}
-
\begin{pmatrix}
 q\\x
\end{pmatrix}
=
g\cdot\left[
 \begin{pmatrix}
  p\\ s x
 \end{pmatrix}
+
\begin{pmatrix}
 q\\x
\end{pmatrix}
\right]\quad\,\, \text{or}\quad \,\,
  \begin{pmatrix}
  p-q\\ (s-1) x
 \end{pmatrix}=g\cdot 
   \begin{pmatrix}
  p+q\\ (s+1) x
 \end{pmatrix}
 .
 $$
 Now we must apply the Cayley transform, i.e., set
$q=y-z$, $p=y+z$.
 It is more convenient to set $q=(y-z)/2$, $p=(y+z)/2$, this does not change a result.
 We get
 $$
\begin{pmatrix}
 z\\(s-1)x
\end{pmatrix}=
g\cdot\begin{pmatrix}
y\\
        (s+1) x
       \end{pmatrix},
 $$
 recall that we consider pairs 
 $z$, $y$, for which there exists  $x$ such that this equality holds.  
set $x'=(s-1)x$. If we change $x$ by $x'$ the condition for  $q$, $p$ does not change.
Thus, 
  $$
\begin{pmatrix}
 z\\x'
\end{pmatrix}=
g\cdot\begin{pmatrix}
y\\
        \frac{s+1}{s-1} x'
       \end{pmatrix}.
 $$
 We get the definition of the characteristic function of the matrix
  $g$ at the point $\frac{s+1}{s-1}$.
 \hfill $\square$

 \sm
 
 {\bf\punct Start of proof of the theorem. The first Cayley transform.} 
 Now, we intend to watch step-by-step pushforwards of the Haar measure
 under the transformations described in Lemma
  \ref{l:2}. It is important that all our manipulations
  over matrices commute with the action of the group $\U(m)$.
 
 The image of the probabilistic Haar measure on
 $\U(k)$ under the Cayley transform was evaluated by Hua Loo Keng.
 It has the form  (see \cite{Hua}, \S 3.1)
 \begin{equation}
 \tau_k \det(1-X^2)^{-k} d\dot X,
 \label{eq:hua}
 \end{equation}
 where $X$ is the Cayley transform of a matrix $g$,
 \begin{equation}
 d\dot X:=\prod_{1\le k\le l \le n} d \Im x_{kl}\cdot \prod \prod_{1\le k< l \le n} d \Re x_{kl},
 \label{eq:dX}
 \end{equation}
 and
 $\tau_k$ is a normalizing constant,
 $$
 \tau_k= 2^{k^2-k}
 \pi^{-k(k+1)/2}\prod_{j=1}^{k-1} j!
 $$
 Obviously, the density of the measure can be written as 
 $$
  \tau_n \det(1+X)^{-k}\det(1-X)^{-k}= \tau_n \bigl|\det(1+X)|^{-2k}.
  $$
  
  In our case,
$k=n+m$. Representing an anti-Hermitian matrix  $H$ (i.e., the Cayley transform
of the matrix  $g$) in a block form
 $H=i\begin{pmatrix}
     A&B\\B^*&D
    \end{pmatrix}$,
    we get a measure
    $$
    \tau_{m+n}\cdot \left|\det\left[1+i\begin{pmatrix}
     A&B\\B^*&D
    \end{pmatrix}\right]\right|^{-2n-2m} \,d\dot A\, d\dot B\, d\dot D,
    $$
 where $d\dot A$, $d\dot B$, $d\dot D$ are the natural Lebesgue measures
 on the space of matrices,
 $$
 d\dot B:=\prod_{1\le k\le n} \prod_{1\le l\le m}\, d\Re b_{kl}\, d\Im b_{kl},
 $$
 and the notations $d\dot A$, $d\dot D$ are similar (\ref{eq:dX}).
 
 \sm
 
    {\bf\punct Quotients with respect to the actions of the groups.%
    \label{ss:quotient}}
   Reduce  $D$ to the diagonal form, let $\mu_1>\dots>\mu_m$
   be the eigenvalues, let
   $M$ be the diagonal form, i.e.,  
   $M=V^{-1} D V$. Denote by $\Xi_m\subset \R^m$ the set of all such collections.
   The distribution of eigenvalues of Hermitian matrices is (see, e.g., \cite{Hua}, \S3.3, or \cite{Far}, Theorem 10.1.4)
\begin{equation}
 dw^m =  \frac{\pi^{m(m-1)/2}}{\prod_{1\le j\le m}j!}  \prod_{1\le k < l\le m} |\mu_k-\mu_l|^2 \prod_{k=1}^m d\mu_k
 .
 \label{eq:vander}
    \end{equation}
    Consider the action of the group 
    $\U(m)$ on the space of Hermitian matrices 
    $\frac 1i H$,
    \begin{equation}
    \begin{pmatrix}
     1&0\\0&V
    \end{pmatrix}
    \begin{pmatrix}
     A&B\\B^*&D
    \end{pmatrix}
  \begin{pmatrix}
     1&0\\0&V
    \end{pmatrix}^{-1}=   \begin{pmatrix}
     A&BV\\(BV)^*&V^{-1}DV
    \end{pmatrix}
    .
    \label{eq:((((}
    \end{equation}
    In fact, this group acts on pairs
    $(D,B)$,
    the matrix $B$ can be regarded as a collection of  
    $n$ rows, denote them by $\beta_1$,\dots, $\beta_n$.
    We get an action of the unitary group
   $\U(m)$ on the collections of the form: a self-adjoint operator
   $D$ and $n$-tuple of vectors 
    $\beta_1$,\dots, $\beta_n$.
   
We can reduce    
$D$ to a diagonal form. After this we preserve a freedom to conjugate
in (\ref{eq:((((})
by elements of the diagonal subgroup  $\T^m$. In this way, we can make all coordinates
of the vector 
$$\beta_1=(b_{11}, \dots, b_{1m})$$
being real and positive. After this all coordinates of vectors $\beta_1$, \dots, $\beta_n$
are fixed.

\begin{lemma}
 Consider the map
 \begin{equation}
\Pi: \Herm_m\times \Mat_{n,m} \to \Xi_m\times (\R_+)^m\times \Mat_{n-1,m},
 \end{equation}
 which to each pair
 $(D,B)$ assigns its canonical form  $(M,B')$.
 Then the image of the Lebesgue measure under the map
 $\Pi$ is
 $$
 dw^m\cdot (2\pi)^m\cdot
  \prod_{k=1}^m b_{1k}'\, db_{1k}' \prod_{1\le l\le n, 2\le k\le n} d\Re b_{kl}'\,d\Im b_{kl}'.
 $$
\end{lemma}

{\sc Proof.} Denote by   $\Fl_m$ the space, whose points are  ordered collections
 $\ell=(\ell_1,\dots,\ell_m)$ of pairwise orthogonal (complex) lines in   $\C^m$. 
 Obviously,  $\Fl_m$ is  simply a flag space). The group $\U(m)$
 acts on this space transitively, a stabilizer of a point is isomorphic
$\T^m$, i.e., $\Fl_n\simeq \U(m)/\T^m$. Since  $\U(m)$ is compact,
there is a unique $\U(m)$-invariant probabilistic measure on
 $\Fl_n$, denote it by  $d\dot \ell$.
 For any Hermitian matrix
 $D$, we assign the collection of its eigenvalues 
$\mu_1> \dots > \mu_m$ and the collection of its eigenvalues  $\ell_1$, \dots, $\ell_m$.
In this way, we get map 
$$
\Herm_m\to \Xi_m\times \Fl_m,
$$
which is defined a.s. and bijective a.s.
The image of the Haar measure under this map is  
 $dw^m\, d\dot\ell$.
This statement is equivalent to%
\footnote{Actually, proofs of formula (\ref{eq:vander})
usually use this fact (directly or hidden).}
formula (\ref{eq:vander}).

Thus, the action of the group
 $\U(m)$ on
 $$\Herm_m\times \Mat_{n,m}\simeq \Herm_m\times (\C^m)^{n}
 $$
 can be regarded as an action on the space
$$
\Xi_m\times\Fl_m\times \C^m\times (\C^m)^{n-1},
$$
trivial on the first factor.

Consider the quotient of
$\Fl_m\times \C^m$ with respect to $\U(m)$. 
Fix a point $\ell\in \Fl_m$, denote by   $\T_\ell^m\subset \U(m)$ 
its stabilizer. For all $\ell$ the distribution of  $\T_\ell^m$-invariants in a fiber $\C^m$ is the same.
Therefore, on the space
$$
(\Fl_m\times \C^m)/\U(m)\simeq (\R_+)^m
$$
we get the same distribution of invariants,
$$
\prod_{k=1}^m b_{1k}'\, db_{1k}'
.
$$

Now we can identify the following spaces by
a $\U(m)$-equivariant measure preserving transformation
  $$
  \Fl_m\times \C^m \simeq (\R_+)^m \times \U(m).
  $$
  Thus we come to the action of
 $\U(m)$ on
$$
\Xi_m\times (\R_+)^m  \times \U(m)\times (\C^m)^{n-1}
,$$
moreover, the action is trivial on the first two factors, on 
$\U(m)$  the group acts by right shifts, on vector-rows it acts in a natural way. 
Now a description of quotient becomes trivial. 
\hfill $\square$

\sm
 
 In fact, we use a slightly modified version of the lemma.
Consider a map
 $$
 \pi:\,\Xi_m\times \Mat_{n,m}\to \Xi_m\times (\R_+)^m\times \Mat_{n-1,m},
 $$
 defined as a reduction of the first row to the canonical form using the action of the torus
$\T^m$ .
 
 \begin{corollary}
 The image of the measure
  $\Xi_m\times \Mat_{n,m} $ under the map  $\pi$ coincides with the image of the measure 
  on $\Herm_m\times \Mat_{n,m}$  under the map 
  $\Pi$. In other words, the natural identification of the spaces 
  $$
  (\Herm_m\times \Mat_{n,m})/\U(m)\longleftrightarrow (\Xi_m\times \Mat_{n,m})/\T^m
  $$
  is measure preserving.
 \end{corollary}

 We came to the following problem.
Consider the space
$\Herm_n\times \Mat_{n,m}\times \Xi_m$
equipped with
 a measure
 \begin{equation}
 \tau_{m+n} \left|\det\left(1+i\begin{pmatrix}
                          A&B\\B^*&M
                         \end{pmatrix}\right)\right|^{-2n-2m}
 dw_m(M) \,d\dot A\,d\dot B
 \label{eq:ABM}
. 
\end{equation}
We must watch behavior of this measure under transformations 2)-4)
described in Lemma
\ref{l:2}.

\sm

{\bf\punct The taking of a characteristic function.%
\label{ss:poles}}
Next, we write the characteristic function of the matrix
$i\begin{pmatrix}
                          A&B\\B^*&M
                         \end{pmatrix}$,
\begin{equation}
 \phi(s)=iA-sB(1-isM)^{-1} B^*=iA-\sum_{k=1}^m \frac{s b_kb_k^*}{1-is\mu_k}
,\end{equation}
where $b_k$  are columns of the matrix  $B$.

We see that poles of $\phi(s)$ are located at points 
 $s_k=\frac 1{i\mu_k}$, and  residues at the poles are rank 1 matrices given by
$$
\Res_{s=1/{i\mu_k}}\phi(s)=-\frac {b_kb_k^*}{\mu_k^2}
.$$
Notice that 
 $\phi(0)=iA$, and the matrix  $b_kb_k^*$ remembers the vector  $b_k$ up to a phase.
Therefore, there is no need to introduce new coordinates,
the matrices $A$, $B$, $M$  can be automatically reconstructed from
the function
$\phi(s)$.

 Notice that in a general position all the vectors
 $b_k$ are non-zero.
 
 \sm
 
 {\bf\punct The second Cayley transform.}
 Next, we take the function
 $$
 \psi(s)=-1+2(1+\phi(s))^{-1}.
 $$
 Notice that for a pure imaginary 
 $s$ the matrix 
 $\psi(s)$ is unitary (since the matrix  $\phi(s)$ is anti-Hermitian).
 Denote
 $U:=\psi(0)$. Clearly, 
 $$
 iA=-1+2(1+U)^{-1}.
 $$
 Next, at points 
 $s=1/i\mu_k$ a matrix   $(1+\phi(s))^{-1}$
 is degenerate, wherefore $(-1)$ is an eigenvalue of $\psi(s)$.
 Let us write  Taylor expansions of  functions 
 $1+\phi(s)$ and $(1+\phi(s))^{-1}$
 at a  point $s=1/i\mu_k$,
 $$
 1+\phi(s)=-\frac {b_kb_k^*}{\mu_k^2\left(s-\frac1{i\mu_k}\right)} + V+\dots,\qquad
 (1+\phi(s))^{-1}=W+\left(s-\frac1{i\mu_k}\right)Y+\dots
 ,$$
 where $V$, $W$, $Y$ are certain matrices. Then
 $$
  (1+\phi(s))^{-1} (1+\phi(s))=1= -\frac {Wb_kb_k^*}{\mu_k^2\left(s-\frac1{i\mu_k}\right)}+ (WV- \mu_k^{-2}Yb_k b_k^*)+...
 $$
 Therefore, 
 \begin{align}
 & Wb_kb_k^*=0;
 \nonumber
 \\
 & WV- \mu_k^{-2}Yb_k b_k^*=1.
 \label{eq:WV}
 \end{align}

 We have 
 $Wb_k=0$. Otherwise, multiplying the vector column  $Wb_k$ by a non-zero row  $b_k$
 we get a nonzero matrix. Next, we write the equality of matrix elements  
$\la (\dots) b_k,b_k\ra$ for both sides  of (\ref{eq:WV}),
\begin{equation}
\bigl\la (WV- \mu_k^{-2}Yb_k b_k^*)b_k, b_k\bigr\ra=\la b_k,b_k\ra.
\label{eq:VWY}
\end{equation}
Notice that  
$$
\la WV b_k,b_k\ra= \la V b_k, W^* b_k\ra=\la V b_k, 0\ra=0.
$$
Indeed, the matrix 
 $W=(\psi(s)+1)/2$ is normal, therefore   $\ker W=\ker W^*$.

Next, 
$$
 (Yb_k b_k^*)b_k=  Yb_k (b_k^*b_k)= \la b_k,b_k\ra\, Yb_k.
$$
Hence  equation
 (\ref{eq:VWY}) transforms to 
$$
\la Y b_k,b_k\ra=-2\mu_k^2,
$$
or, equivalently, to 
\begin{equation}
\la \psi'(1/i\mu_k),b_k\ra=-\mu_k^2.
\label{eq:bb}
\end{equation}
In particular, we see that vectors
 $b_k$ can be reconstructed up to a phase from the function 
 $\psi(s)$.

\sm

{\bf\punct The linear fractional change of the argument.}
It remains to perform the last step, i.e., passing to the function
$$
\chi(t)=\psi\Bigl(\frac{t+1}{t-1}\Bigr).
$$
Let us define new parameters of conjugacy classes:

\sm

1) points $t_k$ of the unit circle defined from the condition:
$$
t_k:=\frac{1+i\mu_k}{1-i\mu_k}\qquad \text{or}\qquad \frac{t_k-1}{t_k+1}=i\mu_k;
$$

2) a matrix  $U\in \U(n)$,
$$
 U:=-1+2(1+iA^{-1})^{-1};
$$

3) vectors $c_k\in \C^n$,
\begin{equation}
 c_k:= \frac12 e^{i\theta_k} (1+t_k)b_k, 
 \label{eq:c-k}
 \end{equation}
 where a phase factor 
 $e^{i\theta_k}$ is fixed by the condition $c_k^1>0$
(where $c_k^1$ is the first coordinate).

Clearly,

\sm

$\bullet$   $U=\chi(-1)$;
   
   \sm
   
$\bullet$   $t_k$ are the points, where  $\chi(t)$ has an eigenvalue $-1$;
   
   \sm
   
$\bullet$   $c_k$ are normalized solutions of the equation $\chi(t_k) v_k=-v_k$.
   
   \sm
   
   Compose a diagonal matrix 
   $T$ with entries $t_j$;
   compose a matrix $C$ of size  $n\times m$ from the columns
 $c_k$.

 The condition (\ref{eq:bb}) can be written as  
 $$
 \la \chi'(t_k) b_k,b_k\ra=-\frac{4}{(t_k+1)^2}=-\frac{4t_k^{-1}}{|t_k+1|^2}
 $$
 or
  $$
 \la \chi'(t_k) c_k,c_k\ra=-t_k^{-1}.
 $$
 
 Now we can watch what is happened with the measure
(\ref{eq:ABM}) under this change.
Simultaneously, we pass to quotient by the action
of $\T^m$, recall that in
 Subsect.\ref{ss:quotient} this was leaved for a future.

 \sm
 
 a) The Lebesgue measure  $d\dot B$ after  change (\ref{eq:c-k}) and passing
 to quotient by the action the torus transforms to
 \begin{equation}
 2^{2m(n-1)} \cdot (2\pi)^{m} \prod_{k=1}^m |1+t_k|^{-2n} 
  \prod_{k=1}^m c_{1k}\, dc_{1k}\, \prod_{1\le l\le n, 2\le k\le n} d\Re c_{kl}'\,\,d\Im c_{kl}'.
  \label{eq:raz}
 \end{equation}
 
b) The measure $dw^m $, see (\ref{eq:vander}), transforms to 
 \begin{equation}
  2^{n(n-1)/2} (2\pi)^m  \cdot 2^{m(m+1)} \prod_{k=1}^m|1+t_k|^{-2m} \prod_{1\le k< l\le m} |t_k-t_l|^2
  \prod_{k=1}^m \frac{dt_k}{i t_k}.
  \label{eq:dva}
  \end{equation}
  
  c) It remains to examine the factor  
   $$
 \tau_{m+n} \left|\det\left(1+i\begin{pmatrix}
                          A&B\\B^*&M
                         \end{pmatrix}\right)\right|^{-2n-2m}
 \,d\dot A
. $$
First, we transform it to the form 
\begin{multline}
 \tau_{m+n} \left|\det\left(1+i\begin{pmatrix}
                          A&B\\B^*&M
                         \end{pmatrix}\right)\right|^{-2n-2m}
 \,d\dot A=\\=
 \frac{ \tau_{m+n}}{\tau_n} \left|\det\left(1+i\begin{pmatrix}
                          A&B\\B^*&M
                         \end{pmatrix}\right)\right|^{-2n-2m} |\det(1+iA)|^{2n}\times\\ \times                          
 \tau_n |\det(1+iA)|^{-2n} d\dot A=\\=
  \frac{ \tau_{m+n}}{\tau_n} \left|\det\left(1+i\begin{pmatrix}
                          A&B\\B^*&M
                         \end{pmatrix}\right)\right|^{-2n-2m} |\det(1+iA)|^{2n}
                         \,d\sigma_n(U)
                         \label{eq:det-long}
, \end{multline}
here we applied the Hua formula 
(\ref{eq:hua}) for the Haar measure.
 
 Next, we must transform two determinants in this formula,
 $$
 \det(1+iA)=\det\bigl(1-1+2(1+U)^{-1}\bigr)=2^n \det(1+U)^{-1}.
 $$
 \begin{multline*}
  \det\left(1+i\begin{pmatrix}
                          A&B\\B^*&M
                         \end{pmatrix}\right)=
                        \det \begin{pmatrix}
                         1 +iA&iB\\iB^*&1+iM
                         \end{pmatrix}
                         =\\=
 \det\begin{pmatrix}
      2(1+U)^{-1}& 2iC(1+T)^{-1}\\
      2i(1+T^{-1})C^*&1+(T+1)^{-1}(T-1)
     \end{pmatrix}                      
 \end{multline*}
 The right lower block equals to
$2T(1+T)^{-1}=2(1+T^{-1})^{-1}$, therefore, we get 
 \begin{multline*}
 2^{m+n} \det \begin{pmatrix}
   1&0\\0&(1+T^{-1})^{-1}
  \end{pmatrix}
 \det\begin{pmatrix}
      (1+U)^{-1}&iC\\
      iC^*& 1+T
     \end{pmatrix}
 \det \begin{pmatrix}
   1&0\\0&(1+T)^{-1}
  \end{pmatrix}
  =\\=
  2^{m+n} \prod_{k=1}^m |1+t_k|^{-2}
  \cdot\det(1+U)^{-1} \det\bigr(1+T+C^*(1+U)C\bigl).
 \end{multline*}
 In the last pass, we applied the formula for a determinant of a block matrix,
 \begin{equation}
 \det\begin{pmatrix}
      a&b\\c&d
     \end{pmatrix}= \det a\, \det(d-ca^{-1} b).
     \label{eq:block-determinant}
 \end{equation}
 Collecting all the formulas together, we come to Theorem
\ref{th:main}.
\hfill $\square$

\sm

Lemma \ref{l:1} does not need a separate proof,
since our step-by-step calculation was based on passing from one parametrization
of double coset space to another.

\noindent
\tt
 Math. Dept., University of Vienna; \\
ITEP, Moscow \\
MechMath, Moscow State University\\
IITP, Moscow
\\
neretin(at)mccme.ru
\\
URL:http://www.mat.univie.ac.at/$\sim$neretin/

\vspace{60pt}

Key words: inner functions, characteristic functions, Haar measure, Cayley transform, random functions

UDC 512.546.32, 517.547.5, 517.548.5

MSC 47A48, 28C10, 20E45
\end{document}